\documentclass[oneside]{amsart}

\theoremstyle{plain}    
\newtheorem{thm}{Theorem}[section]
\numberwithin{equation}{section} 
\numberwithin{figure}{section} 
\theoremstyle{plain}    
\newtheorem{cor}[thm]{Corollary} 
\theoremstyle{plain}    
\newtheorem{lem}[thm]{Lemma} 
\theoremstyle{plain}    
\newtheorem{prop}[thm]{Proposition} 
\theoremstyle{remark}
\newtheorem{rem}[thm]{Remark}

\begin{document}

\title[Mirror symmetry for projective spaces]{Semi-infinite Hodge structures and mirror symmetry for projective spaces}\thanks{Date: June 27, 2000; revised (a few local improvements and references added) September 6, 2001; Preprint ENS DMA-00-24}

\author{Serguei Barannikov}

\address{\'{E}cole normale sup\'{e}riere (Paris), UMR8553 CNRS\\
 }

\email{serguei.barannikov@ens.fr }

\urladdr{http//www.dma.ens.fr/\protect\( \sim \protect \)barannik}

\begin{abstract}
We express total set of rational Gromov-Witten invariants of \( \mathbb {C}\mathbb {P}^{n} \)
via periods of variations of semi-infinite Hodge structure associated with their
mirror partners. For this explicit example we give detailed description of general
construction of solutions to WDVV-equation from variations of semi-infinite
Hodge structures of Calabi-Yau type which was suggested in a proposition from
our previous paper (\cite{B2} proposition 6.5). 
\end{abstract}
\maketitle
\setcounter{tocdepth}{1}

\tableofcontents{}

\section{Introduction}

Motivated by higher-dimensional mirror symmetry we introduced in \cite{B2}
quantum periods associated with non-commutative deformations of complex manifolds.
When it is just a deformation of complex structure then the corresponding quantum
periods specialize to usual periods. In another simplifying situation quantum
periods reduce to oscillating integrals. In this paper we look closely at this
case and associated variations of semi-infinite Hodge structures. It turns out
that these variations give mirror partners to projective spaces. Our principal
result is the theorem \ref{th:main} expressing the generating function for
total collection of rational Gromov-Witten invariants of \( \mathbb {C}\mathbb {P}^{n} \)
via periods of variations of semi-infinite Hodge structure of its mirror partner.

Let us consider a pair \( (X,f) \) where 
\[
X=\{x_{0}\cdot x_{1}\cdot \ldots \cdot x_{n}=1\}\subset \mathbb {C}^{n+1},\, \, f:X\rightarrow \mathbb {A}_{\mathbb {C}}^{1},\, \, f=x_{0}+\ldots +x_{n}\]
 Such data was conjectured (see \cite{G}, \cite{EHX}) to be mirror partner
of projective space \( \mathbb {C}\mathbb {P}^{n} \). Generating function (potential)
for \( genus=0 \) Gromov-Witten invariants of \( \mathbb {C}\mathbb {P}^{n} \)
reads as (see \cite{KM} eq.(5.15)) 
\begin{equation}
\mathcal{F}^{\mathbb {C}\mathbb {P}^{n}}=\frac{1}{6}\sum _{0\leq i,j,k\leq n}y^{i}y^{j}y^{k}\delta _{i+j,n-k}+\sum _{d;m_{2},\ldots ,m_{n}}N(d;m_{2},\ldots ,m_{n})\frac{(y^{2})^{m_{2}}\ldots (y^{n})^{m_{n}}}{m_{2}!\ldots m_{n}!}e^{dy^{1}}
\end{equation}
 here \( y\in H^{*}(\mathbb {C}\mathbb {P}^{n},\mathbb {C}) \) and the Taylor
coefficients \( N(d;m_{2},\ldots ,m_{n}) \) are the virtual numbers of rational
curves of degree \( d \) in \( \mathbb {C}\mathbb {P}^{n} \) intersecting
\( m_{a} \) hyperplanes of codimension \( a \).

Let us put 
\begin{equation}
\label{hatF}
\widehat{F}(x;t)=f+\sum ^{n}_{m=0}t^{m}(\sum ^{n}_{i=0}x_{i})^{m}
\end{equation}
 An oscillating integral \( \psi _{k}=\int _{\Delta _{k}}\text {exp}\, (\frac{\widehat{F}(x;t)}{\hbar })(a_{0}+\hbar a_{1}+\ldots )\frac{dx_{0\ldots }dx_{n}}{d(x_{0\ldots }x_{n})} \)
can be expressed using integration by parts via integrals 
\[
\varphi ^{m}_{k}(t,\hbar )=\int _{\Delta _{k}}\text {exp}\, (\frac{\widehat{F}(x;t)}{\hbar })(\sum ^{n}_{i=0}x_{i})^{m}\frac{dx_{0\ldots }dx_{n}}{d(x_{0\ldots }x_{n})},\, \, \, \varphi ^{m}_{k}=\hbar \frac{\partial \varphi ^{0}_{k}(t,\hbar )}{\partial t^{m}}\]
 as linear combination \( \psi _{k}=\sum _{m=0}^{n}v_{m}\varphi ^{m}_{k} \)
with coefficients \( v_{m}=v_{m}^{(0)}+\hbar v_{m}^{(1)}+\ldots  \). Let us
consider an element \( \psi (t,\hbar ) \) satisfying the following normalization
condition: 
\[
\psi _{k}(t,\hbar )=\sum _{m=0}^{n}u_{m}\varphi ^{m}_{k}(0,\hbar )\]
 for some functions \( u_{m}(t,\hbar )=\delta _{m,0}+\hbar ^{-1}u_{m}^{(-1)}(t)+\ldots  \),
\( m=0,\ldots ,n \) (exact meaning of this condition is given in \S \ref{sect:WDVV}).
We will show that if one makes a change of parameters \( y^{m}=u^{(-1)}_{m}(t) \)
then Picard-Fuchs equations for oscillating integrals \( \psi _{k} \) take
the form 
\begin{equation}
\label{Aabcprim}
\frac{\partial ^{2}\psi _{k}}{\partial y^{i}\partial y^{j}}=\hbar ^{-1}\sum _{k}A_{ij}^{m}(y)\frac{\partial \psi _{k}}{\partial y^{m}}
\end{equation}
 Our principal formula can be written then as 
\begin{equation}
\forall \, i,j,m,\, \, \, A_{ij}^{(n-m)}(y)=\partial ^{3}_{ijm}\mathcal{F}^{\mathbb {C}\mathbb {P}^{n}}(y)
\end{equation}
in other words we give integral representation for structure constants of \emph{full}
quantum cohomology algebra of \( \mathbb {C}\mathbb {P}^{n} \).

In fact, the functions \( \psi _{k}(t,\hbar ) \) should be viewed as periods
of variation of semi-infinite Hodge structures attached to \( (X,f) \). The
semi-infinite analog of Hodge filtration associated with a point \( t \) of
moduli space of deformations of \( f \) is given by subspace generated by cohomology
classes of oscillating integrals in the space of sections of local system over
\( \mathbb {C}\setminus \{0\} \) whose fiber over \( \hbar  \) is the relative
cohomology group \( H^{n}(X,\text {Re}\, \frac{\widehat{F}(x;t)}{\hbar }\rightarrow -\infty ;\mathbb {C}) \).
The normalization condition imposed on element \( \psi (t) \) from this subspace
means that it should belong at the same time to certain covariantly constant
affine subspace.

Another aim of this paper is to consider in this explicit case the general relation,
which we described in (\cite{B2}, proposition 6.5) between variations of semi-infinite
Hodge structures and Frobenius manifolds. We construct in particular in \S\ref{sect:WDVV}
a family of solutions to WDVV equations parameterized by open part in isotropic
affine grassmanian. Same result holds true for similar pairs \( (X,f) \) satisfying
homological condition of Calabi-Yau type, see remark after proposition \ref{pr:frob}.
The approach to WDVV-equations described in this paper arose from an attempt
to understand relation of the construction of solutions to these equations from
(\cite{B1}, section 5) with elegant geometric approach to integrable hierarchies
from \cite{SW}.

Concerning the previously known results we should mention the paper \cite{M}
dealing with the case of \( \mathbb {C}\mathbb {P}^{2} \) in which after some
change of variables the associativity equation \( \widetilde{\mathcal{F}}_{y^{2}y^{2}y^{2}}=\widetilde{\mathcal{F}}^{2}_{y^{1}y^{1}y^{2}}-\widetilde{\mathcal{F}}_{y^{1}y^{1}y^{1}}\widetilde{\mathcal{F}}_{y^{1}y^{2}y^{2}} \)
for non-trivial piece of the potential \( \mathcal{F}^{\mathbb {C}\mathbb {P}^{2}} \)
was identified with differential equation on (multi)sections of pencil of elliptic
curves \( z^{2}=x(x-1)(x-t) \) written in terms of \( \int _{\infty }^{x(t)}\frac{dx}{z} \).
Other relevant important result is an integral representation given in (\cite{G})
for structure constants of \emph{small} quantum cohomology algebra \( \mathbb {C}[p,q]/p^{n+1}=q \)
of \( \mathbb {C}\mathbb {P}^{n} \) and also for its \( T^{n}- \)equivariant
counterpart. In terms of enumerative geometry of curves in \( \mathbb {C}\mathbb {P}^{n} \)
Taylor expansion of structure constants of this algebra encode the simplest
number \( N(1;0,\ldots ,0,2)=1 \) which is the number of lines going through
two points in \( \mathbb {C}\mathbb {P}^{n} \).

\section{Moduli space.}

In this section we consider meaning of parameter space of deformation (\ref{hatF})
from the point of view of deformation theory. The content of this section is
not used below and readers interested only in explicit formulas giving integral
representations for Gromov-Witten invariants may skip it.

Deformations of the pair \( (X,f) \) can be described in general with help
of differential graded Lie algebra of polyvector fields 
\[
\frak {g}(X,f)=\oplus _{i}\frak {g}^{i}(X,f)[-i],\, \, \, \frak {g}^{i}(X,f):=\Gamma _{Zar}(X,\Lambda ^{1-i}T),\, \, \, d:=[f,\cdot ]\]
 equipped with the Schouten bracket. Explicitly, the elements \( \gamma \in \frak {g}(X,f) \)
can be uniquely written in the form 
\[
\gamma =\sum _{I\subseteq \{1,\ldots ,n\}}\gamma _{I}(x_{1},x^{-1}_{1},\ldots ,x_{n},x^{-1}_{n})\partial _{i_{1}}\wedge \ldots \wedge \partial _{i_{k}}\]
 where \( \gamma _{I}\in \mathbb {C}[x_{1},x^{-1}_{1},\ldots ,x_{n},x^{-1}_{n}] \).
The space of first order deformations at the base point \( T_{[(X,f)]}\mathcal{M}^{(X,F)} \)
equals to total cohomology group of complex \( (\frak {g}(X,f),[f,\cdot ]) \).
Simple calculation shows that the cohomology of complex \( (\frak {g}(X,f),[f,\cdot ]) \)
are zero everywhere except in degree zero where they are \( (n+1)- \)dimensional.
Notice that the partial derivatives \( \frac{\partial \widehat{F}}{\partial t_{i}}|_{t=0} \)
\( i=0\ldots n \), where \( \widehat{F}(x;t_{0},\ldots ,t_{n}) \) is the family
of functions defined in (\ref{hatF}), form a basis in the cohomology of the
complex \( (\frak {g}(X,f),[f,\cdot ]) \). It follows that \( \frak {g}(X,f) \)
is formal as differential graded Lie algebra and that the set \( \mathcal{M}^{(X,F)}(R) \)
of equivalence classes of solutions to Maurer-Cartan equation in \( \mathfrak {g}^{1}\otimes R \)
for an Artin algebra \( R \) can be identified with 
\[
\mathcal{M}^{(X,F)}(R)=\{F=f+\widetilde{f},\widetilde{f}\in \Gamma
_{Zar}(X,\mathcal{O}_{X})\otimes \mathfrak {M}_{R}\}/\{\varphi \in
Aut_{Zar}(X)_{/R}, \varphi(0)=Id\}\]
 in other words \( \mathcal{M}^{(X,F)} \) is the moduli space of deformations
of the algebraic function \( f \). It follows that \( \widehat{F}(x;t_{0},\ldots ,t_{n}) \)
is a mini-versal family for such deformations. In the sequel it will be convenient
for us to denote by \( \mathcal{U} \) the germ of \( (n+1)- \)dimensional
smooth analytic space of the base of the deformation \( \widehat{F}(x;t_{0},\ldots ,t_{n}) \)
and by \( \widetilde{X}=X\times \mathcal{U} \) the total space on which \( \widehat{F}(x;t_{0},\ldots ,t_{n}) \)
is defined.

\section{Variation of semi-infinite Hodge structures associated with \protect\( (X,\widehat{F})\protect \).}

In this section we define semi-infinite analog of variation of Hodge structures
associated with the moduli space \( \mathcal{M}^{(X,F)} \). The analogy with
usual variations of Hodge structures associated with families of complex manifolds
is explained in \cite{B2} \S 3.

Let \( \mathcal{R}_{f} \) denotes the sheaf of relative cohomology whose fiber
over \( \hbar \in \mathbb {C}\setminus {\{0\}} \) equals to 
\begin{equation}
\label{Hfh}
\mathcal{R}_{f,\hbar }:=H^{n}(X,\text {Re}\, \frac{f}{\hbar }\rightarrow -\infty ;\mathbb {C})
\end{equation}
 Since there is a natural integral structure on this sheaf: \( \mathcal{R}_{f,\hbar }=H^{n}(X,\text {Re}\, \frac{f}{\hbar }\rightarrow -\infty ;\mathbb {Z})\otimes \mathbb {C} \),
therefore it is equipped with Gauss-Manin connection which we will denote by
\( \mathcal{D}_{\partial /\partial \hbar } \). Denote also via \( (\mathcal{R}_{\widehat{F}(t)},\mathcal{D}_{\partial /\partial \hbar }^{\mathcal{R}_{\widehat{F}(t)}}) \)
the analogous sheaf with connection associated with the function \( \widehat{F}(t) \),
\( t\in \mathcal{U} \). There is a flat connection \( \nabla ^{\mathcal{M}} \),\break  \( [\nabla ^{\mathcal{M}} \),
\( \mathcal{D}_{\partial /\partial \hbar }^{\mathcal{R}_{\widehat{F}(t)}}]=0 \)
on the total family of sheaves with connections \( (\mathcal{R}_{\widehat{F}(t)},\mathcal{D}_{\partial /\partial \hbar }^{\mathcal{R}_{\widehat{F}(t)}}) \),
\( t\in \mathcal{M} \). 

Let us consider subspace of sections of sheaf \( \mathcal{R}_{\widehat{F}(t)} \)
which are represented by oscillating integrals. In other words we consider the
subspace which consists of sections given by integrals of elements of the form
\( [\text {exp}\, (\frac{1}{\hbar }\widehat{F}(t))\sum _{i\geq 0}\hbar ^{i}\varphi _{i}] \),
where \( \sum _{i\geq 0}\hbar ^{i}\varphi _{i}\in \Gamma (X\times \mathbb {A}^{1}_{\mathbb {C},\text {an}},\Omega ^{n}_{X\times \mathbb {A}^{1}_{\mathbb {C},\text {an}}/\mathbb {A}^{1}_{\mathbb {C},\text {an}}}) \)
(here \( \mathbb {A}^{1}_{\mathbb {C},\text {an}} \) denotes the spectrum of
the algebra of analytic functions in \( \hbar  \)). 

A natural basis in \( H_{n}(X,\text {Re}\, \frac{f}{\hbar }\rightarrow -\infty ;\mathbb {C}) \)
can be constructed by means of Morse theory. Namely given a critical point \( p \)
of \( f \) and a choice of metric on \( X \) the set of all gradient lines
\( \rho (s) \) of \( \text {Re}\, \frac{f}{\hbar } \) such that \( \rho (s)\rightarrow p \)
as \( s\rightarrow +\infty  \) form relative cycle \( \Delta ^{+}_{p}(\hbar )\in  \)
\( H_{n}(X,\text {Re}\, \frac{f}{\hbar }\rightarrow -\infty ;\mathbb {C}) \).
One can see easily using standard methods of toric geometry that there exists
partial smooth compactification \( \overline{X} \) of \( X \) such that extension
of \( f \) on \( \overline{X} \) is proper and such that its restricton on
\( \overline{X}\setminus X \) has no critical points. Then it is easy to see
that there exists a metric on \( \overline{X} \) such that gradient flow of
\( \text {Re}\, (f/\hbar ) \) preserves \( \overline{X}\setminus X \). Let
\( \Delta ^{+}_{p_{i}}(\hbar ) \) are relative cycles constructed with help
of restriction of this metric on \( X \), where \( \{p_{i}\} \), \( p_{i} \):
\( x_{0}=\ldots =x_{n}=\exp \frac{2\pi \sqrt{-1}i}{n+1} \), \( i=0,\ldots ,n \)
is the set of critical points of \( f \). Standard arguments of Morse theory
show that \( \{\Delta ^{+}_{p_{i}}(\hbar )\} \) is a basis in \( H_{n}(X,\text {Re}\, \frac{f}{\hbar }\rightarrow -\infty ;\mathbb {C}) \). 

Let us describe explicitly the subspace of oscillating integrals at \( t=0 \).
It is convenient to introduce an auxiliary variable \( \alpha  \). Let us denote
via \( \xi _{k}(\hbar ) \) the coefficient in front of \( \alpha ^{k} \) in
the expression 
\begin{equation}
\label{xi}
\xi (\alpha ,\hbar ):=\exp (-\alpha (n+1)\log \hbar )\sum ^{\infty }_{d=0}\frac{1}{\hbar ^{(n+1)d}[\Gamma (\alpha +d+1)]^{n+1}}
\end{equation}
 where \( \frac{1}{\Gamma (\alpha +d+1)}=\frac{1}{\Gamma (\alpha +1)}\prod ^{d}_{i=1}\frac{1}{i}(1-\frac{\alpha }{i}+(\frac{\alpha }{i})^{2}-\ldots ) \)

\begin{prop}
\label{prL0}
\begin{equation}
\label{L0}
\int _{\Delta _{k}\subset X}\text {exp}\, (\frac{\sum _{i=0}^{n}x_{i}}{\hbar })\frac{dx_{0\ldots }dx_{n}}{d(x_{0\ldots }x_{n})}=\xi _{k}(\hbar ),\, \, \, k=0,\ldots ,n
\end{equation}
 for some locally constant basis \( \{\Delta _{k}(\hbar )\} \) in \( H_{n}(X,\text {Re}\, \frac{f}{\hbar }\rightarrow -\infty ;\mathbb {C}) \). 
\end{prop}
\begin{proof}
Recall the integral representation: 
\[
\frac{1}{\Gamma (s+1)}=\frac{1}{2\pi \sqrt{-1}}\int _{C}z^{-s}\exp z\frac{dz}{z}\]
 where \( C \) is the following contour on the complex plane \( \{z|Im\, z=0-,\, Re\, z\in ]-\infty ,-a]\} \)\( \sqcup \{z||z|=a\}\sqcup  \)\( \{z|Im\, z=0+,\, Re\, z\in ]-\infty ,-a]\} \).\( a>0 \).
Therefore
\[
\xi =\frac{1}{(2\pi \sqrt{-1})^{n+1}}\sum ^{\infty }_{d=0}\int \ldots \int _{C\times \ldots \times C}(\hbar z_{0}\cdot \ldots \cdot \hbar z_{n})^{-d-\alpha }\exp (z_{0}+\ldots +z_{n})\frac{dz_{0\ldots }dz_{n}}{z_{0\ldots }z_{n}}\]
 Notice that because \( \frac{1}{\Gamma (d+1)}=0 \) for \( d<0 \) the first
\( (n+1) \) Taylor coefficients of \( \xi (\alpha ,\hbar ) \) do not change
if one takes \( d=-\infty  \) as lower limit in summation. Now after making
a change of variables \( x_{1}=\hbar z_{1} \), ...,\( x_{n}=\hbar z_{n} \),
\( y=\prod ^{n}_{i=0}(\hbar z_{i}) \) and noticing that \( \sum _{d}y^{-d}=\delta (y-1) \)
in appropriate distributional sense one arrives at \ref{L0}. We leave details
to the interested reader. Similar formulas can be found in \cite{Kapr}. Alternatively
one can notice that both \( \xi _{k}(\hbar ) \) and \( \int _{\Delta ^{+}_{p_{i}}}e^{f/\hbar }dx_{0}\ldots dx_{n}/d(x_{0}\ldots x_{n}) \)
satisfy \( (n+1)- \)st order differential equation \( (-\frac{\hbar }{n+1}\frac{\partial }{\partial \hbar })^{n+1}s=\hbar ^{-(n+1)}s \).
It follows from expansions
\[
\int _{\Delta ^{+}_{p_{i}}}e^{\frac{f}{\hbar }}\frac{dx_{0}\ldots dx_{n}}{d(x_{0}\ldots x_{n})}=e^{\frac{f(p_{i})}{\hbar }}(\text {const}\, \hbar ^{\frac{n}{2}}+O(\hbar ^{\frac{n}{2}+1}))\, \text {as}\, \hbar \rightarrow 0\]
 and \( \xi _{k}(\hbar )=(-(n+1)\log \hbar )^{k}+O((\log \hbar )^{k-1}) \)
as \( \hbar \rightarrow \infty  \) that we have two \( (n+1) \) -tuples of
linearly independent solutions, which therefore must be related by a linear
transformation. 
\end{proof}
Consequently the subspace of sections of sheaf \( \mathcal{R}_{f} \) represented
by oscillating integrals can be described explicitly as subspace generated by
elements which in the locally constant frame dual to \( \{\Delta _{k}\} \)
are written as 
\[
\hbar ^{2l+m}(\frac{\partial ^{l}\xi _{k}(\hbar )}{\partial \hbar ^{l}})_{k=0,\ldots ,n}\, ,\, \, l\in \{0,\ldots ,n\},\, \, m\in \mathbb {N}\cup \{0\}\]
 It is convenient to introduce the algebra \( \mathbb {C}[\alpha ]/\alpha ^{n+1}\mathbb {C}[\alpha ] \)
and to identify \( \{\alpha ^{k}\}_{k=0,\ldots ,n} \) with locally constant
frame dual to \( \{\Delta _{k}\} \).Then monodromy transformation around \( \hbar =0 \)
acting on cohomology \( H^{n}(X,\text {Re}\, \frac{f}{\hbar }\rightarrow -\infty ;\mathbb {C}) \)
is given in the basis \( \{\alpha ^{k}\} \) by multiplication by \( \exp ((n+1)2\pi \sqrt{-1}\alpha ) \).
Notice that the set \( \{\hbar ^{-\alpha (n+1)}\alpha ^{k}\} \) is a single-valued
frame in \( \mathcal{R}_{f} \). Below \( \alpha  \) will always denote an
element satisfying \( \alpha ^{n+1}=0 \).

Let us consider a completion of the space \( \Gamma (\mathbb {A}^{1}_{\mathbb {C},\text {an}}\setminus \{0\},\mathcal{R}_{f}) \)
defined as Hilbert space \( H=\Gamma _{L^{2}}(S^{1},\mathcal{R}_{f}|_{S^{1}}) \)
consisting of \( L^{2}- \)sections of restriction of the sheaf \( \mathcal{R}_{f} \)
to the circle \( S^{1}=\{\hbar |\hbar \in \mathbb {C},\, |\hbar |=R,\, R>0\} \).
Explicitly, elements of \( H \) can be written as 
\[
\hbar ^{-\alpha (n+1)}\sum _{i=0}^{n}\alpha ^{i}\zeta _{i}(\hbar ),\, \, \text {where}\, \, \zeta _{i}(\hbar )\in L^{2}(S^{1},\mathbb {C})\]
 where \( \hbar ^{-\alpha (n+1)}:=\exp (-\alpha (n+1)\log \hbar ) \). Let us
consider Segal-Wilson grassmanian \( Gr \) associated with \( H \). We introduce
polarization \( H=H^{+}\oplus H^{-} \) where \( H^{+} \) and \( H^{-} \)
are the closed subspaces generated by elements of the form \( \hbar ^{-\alpha (n+1)+k}\alpha ^{i} \)
for \( k\geq 0 \) and \( k<0 \) respectively. Recall (see \cite{PS}) that
Segal-Wilson grassmanian consists of closed subspaces \( L\subset H \) \char`\"{}comparable\char`\"{}
with \( H^{+} \): \( pr_{+}:L\rightarrow H^{+} \) is a Fredholm operator and
\( pr_{-}:L\rightarrow H^{-} \) is a Hilbert-Schmidt operator.

The semi-infinite subspace \( L(t) \) associated with a point \( t\in \mathcal{U} \)
is defined as closure of the subspace of oscillating integrals, i.e. it is the
closed subspace generated by elements , which in the frame dual to \( \{\Delta _{k}\} \)
have coordinates 
\begin{equation}
\label{eq:int}
\int _{\Delta _{k}(t,\hbar )}\exp (\frac{\widehat{F}(x;t)}{\hbar })\, \varphi \, \hbar ^{i},\, \, \, i\geq 0,\, \, \varphi \in \Gamma _{Zar}(X,\Omega _{X}^{n+1})
\end{equation}
 where \( \Delta _{k}(t,\hbar ) \), is  covariantly constant family
of elements in \( H^{n}(X,\text {Re}\, \frac{\widehat{F}(x;t)}{\hbar }\rightarrow -\infty ;\mathbb {C}) \)
such that \( \Delta _{k}(0,\hbar ) \) coincides with \( \Delta _{k}(\hbar ) \)
from prop.\ref{prL0}. More precisely we consider \( (n+1)\times (n+1) \) matrix
\( \vartheta  \) whose entries are series of hypergeometric type obtained via
expansion \( \exp (\sum ^{n}_{m=0}t^{m}f^{m}/\hbar )=\sum _{l}s^{l}(t,\hbar )f^{l} \)
and termwise integration:
\[
\vartheta _{kj}=\int _{\Delta _{k}(t)}\exp (\frac{f+\sum ^{n}_{m=0}t^{m}f^{m}}{\hbar })f^{j}\frac{dx_{0}\ldots dx_{n}}{d(x_{0}\ldots x_{n})}=\sum _{l=j}^{\infty }s^{l-j}(t,\hbar )\int _{\Delta _{k}}e^{\frac{f}{\hbar }}f^{l}\frac{dx_{0}\ldots dx_{n}}{d(x_{0}\ldots x_{n})}\]
 where \( \int _{\Delta _{k}}e^{\frac{f}{\hbar }}f^{l}\frac{dx_{0}\ldots dx_{n}}{d(x_{0}\ldots x_{n})}=(\frac{\partial }{\partial (\hbar ^{-1})})^{l}\xi _{k} \)
and \( \xi _{k}(\hbar ) \) is defined in (\ref{xi}). The resulting series
is convergent for \( |\hbar |>C \) for some constant \( C>0 \) and 
small \( t \). We set \( L(t):=\vartheta H^{+} \). It follows from \cite{PS}
prop. 6.3.1 that \( L(t)\in Gr \).

The family of subspaces \( L(t) \) has the following basic property.

\begin{prop}
The family of subspaces \( L(t),t\in \mathcal{U} \) satisfies semi-infinite
analog of Griffiths transversality condition with respect to the Gauss-Manin
connection \( \nabla ^{\mathcal{M}} \): 
\begin{equation}
\label{griftransv}
\forall \, i\, \, \, \frac{\partial }{\partial t^{i}}L(t)\subseteq \hbar ^{-1}L(t)
\end{equation}

\end{prop}
\begin{proof}
It follows from \( \frac{\partial }{\partial t^{i}}\text {exp}(\frac{\widehat{F}(t)}{\hbar })=\frac{1}{\hbar }f^{i}\text {exp}(\frac{\widehat{F}(t)}{\hbar }) \)
\end{proof}
The induced map (symbol of the Gauss-Manin connection): 
\[
Symbol(\nabla ^{\mathcal{M}}):L([f])/(\hbar L([f]))\otimes T_{[f]}\mathcal{M}\rightarrow (\hbar ^{-1}L([f]))/L([f])\]
 is easy to describe. At \( [f]\in \mathcal{M} \) space \(
 L([f])/(\hbar L([f]))\simeq (\hbar ^{-1}L([f]))/L([f])\)
is identified naturally with 
\[
\Gamma _{Zar}(X,\Omega _{X}^{n})/\{d_{X}f\wedge \nu \, |\, \nu \in \Gamma _{Zar}(X,\Omega _{X}^{n-1})\}\]
 and
\begin{equation}
\label{TtM}
T_{[f]}\mathcal{M}=\Gamma _{Zar}(X,\mathcal{O}_{X})/\{Lie_{v}f)\, |\, v\in \Gamma _{Zar}(X,\mathcal{T}_{X})\}
\end{equation}
 Then 
\[
Symbol(\nabla ^{\mathcal{M}})(\varphi \otimes \mu )=[\varphi \cdot \mu ]\]

\section{Solutions to WDVV-equation.}

\label{sect:WDVV}

In this section we will show how to construct a solution to WDVV-equations associated
with \( (X,f) \) given an arbitrary element from an open domain in isotropic
affine grassmanian of \( H \).

In the next section we single out a solution which coincides with solution defined
by Gromov-Witten invariants of \( \mathbb {C}\mathbb {P}^{n} \).

\subsection{Picard-Fuchs equation for periods of the semi-infinite Hodge structure}

We start from fixing a choice of semi-infinite subspace \( S\subset H \),
\( \hbar ^{-1}S\subset S \) transversal to \( L(0) \): 
\begin{equation}
\label{transv}
H=L(0)\oplus S
\end{equation}
\( S \) belongs to the same type of grassmanian related with polarization where
roles of \( H^{+} \) and \( H^{-} \) are inversed.

Next we would like to choose an element \( \Omega _{0}\in L(0) \) such that
the symbol of Gauss-Manin connection restricted to the class \( [\Omega _{0}]\, mod\, (\hbar L([f])) \)
in the quotient \( L([f])/(\hbar \, L([f])) \) gives an isomorphism 
\begin{equation}
\label{Omega-prty}
Symbol(\nabla ^{\mathcal{M}})([\Omega _{0}]\otimes \cdot ):T_{[f]}\mathcal{M}\rightarrow (\hbar ^{-1}L([f]))/L([f])
\end{equation}
 It is easy to see that 
\begin{eqnarray}
T_{[f]}\mathcal{M}\simeq \mathbb {C}[p]/(p^{n+1}-1)\, \text {where}\, p=[(\frac{1}{n+1}\sum ^{n}_{i=0}x_{i})]\notag  &  & \label{Tf} \\
L([f])/(\hbar L([f]))\simeq \left( \mathbb {C}[p]/(p^{n+1}-1)\right) [\frac{dx_{0\ldots }dx_{n}}{d(x_{0\ldots }x_{n})}]\notag  & 
\end{eqnarray}
 So one can take here as \( \Omega _{0} \) for example an arbitrary element
of the form 
\[
\left[ \exp (\frac{1}{\hbar }f)\left( (\sum ^{n}_{i=0}x_{i})^{k}+\hbar \omega _{1}+\ldots \right) \frac{dx_{0\ldots }dx_{n}}{d(x_{0\ldots }x_{n})}\right] \, k=0,\ldots ,n\]
 The transversality condition (\ref{transv}) implies that the intersection
of \( L(t) \) for \( t\in \mathcal{U} \) with affine space \( \{\Omega _{0}+\eta |\eta \in S\} \)
consists of a single element. Let us denote it via \( \Psi ^{S} \): 
\[
\{\Psi ^{S}(t,\hbar )\}=L(t)\cap \{\Omega _{0}+\eta |\eta \in S\}\]
 Notice that the property (\ref{Omega-prty}) and the mini-versality of the
family \( \widehat{F}(t) \) implies that \( \{\partial _{i}\Psi ^{S}(t)\in \hbar ^{-1}L(t)\, \text {mod}\, L(t)\} \)
is a basis in $(n+1)-$dimensional vector space \( (\hbar ^{-1}L(t))/L(t) \)
and therefore \( \{\partial _{i}\Psi ^{S}(t)\in S\, \text {mod}\, \hbar ^{-1}S\} \)
is a basis in \( S/(\hbar ^{-1}S) \). Hence
the map 
\begin{equation}
\label{Psi(infty)}
\Psi ^{S}(t,\infty ):=[\Psi ^{S}(t,\hbar )-\Omega _{0}]\in S\, \text {mod}\, \hbar ^{_{-1}}S
\end{equation}
 is a local isomorphism and induces a set of coordinates on \( \mathcal{U} \)
which we denote by \( \{t_{S}^{a}\} \). Let us denote via \( \Psi ^{S,k}=\int _{\Delta _{k}(t_{S}^{a},\hbar )}\Psi ^{S} \)
the components of \( \Psi ^{S} \) with respect to the basis \( {\alpha ^{k}} \).

\begin{prop}
\label{Pr_pic-fuchs} The periods \( \Psi ^{S,k}(t_{S}^{a},\hbar ) \) satisfy
\begin{equation}
\label{Aabc}
\frac{\partial ^{2}\Psi ^{S,k}}{\partial t_{S}^{a}\partial t^{b}_{S}}=\hbar ^{-1}\sum _{c}A_{ab}^{c}(t_{S})\frac{\partial \Psi ^{S,k}}{\partial t_{S}^{c}}
\end{equation}

\end{prop}
\begin{proof}
It follows from the prop. \ref{griftransv} that 
\[
\frac{\partial ^{2}\Psi ^{S}}{\partial t^{i}\partial t^{j}}\in \hbar ^{-2}L(t)\]
 Therefore one has 
\[
\frac{\partial ^{2}\Psi ^{S}}{\partial t^{i}\partial t^{j}}-\hbar ^{-1}\sum _{m}\left( A^{(-1)}\right) _{ij}^{m}(t)\frac{\partial \Psi ^{S}}{\partial t^{m}}-\sum _{m}\left( A^{(0)}\right) _{ij}^{m}(t)\frac{\partial \Psi ^{S}}{\partial t^{m}}\in L(t)\]
 for some \( \left( A^{(-1)}\right) _{ij}^{m}(t),\, \left( A^{(0)}\right) _{ij}^{m}(t)\in \mathcal{O}^{analytic}_{\mathbb {C}} \).
On the other hand, 
\[
\frac{\partial ^{2}\Psi ^{S}}{\partial t^{i}\partial t^{j}},\frac{\partial \Psi ^{S}}{\partial t^{i}}\in S\]
 Therefore 
\[
\frac{\partial ^{2}\Psi ^{S}}{\partial t^{i}\partial t^{j}}=\hbar ^{-1}\sum _{m}\left( A^{(-1)}\right) _{ij}^{m}(t)\frac{\partial \Psi ^{S}}{\partial t^{m}}+\sum _{m}\left( A^{(0)}\right) _{ij}^{m}(t)\frac{\partial \Psi ^{S}}{\partial t^{m}}\]
 The coordinates \( \{t^{a}_{S}\} \) were chosen so that \( \Psi ^{S}\, \text {mod}\, \hbar ^{-1}S:\mathcal{U}\rightarrow S/\hbar ^{-1}S \)
is linear in \( \{t^{a}_{S}\} \). Therefore in these coordinates \( \left( A^{(0)}\right) ^{c}_{ab}=0 \). 
\end{proof}
\begin{cor}
\label{cor:AA,dA} One has \( [A,A]=0 \), \( dA=0 \) for \( A=\sum _{a}A_{ab}^{c}(t_{S})dt^{a}_{S} \). 
\end{cor}
In particular, the following formula defines commutative and associative product
\char`\"{}\( \circ  \)\char`\"{} on tangent sheaf of \( \mathcal{U} \):
\begin{equation}
\label{circ}
\frac{\partial }{\partial t^{a}}\circ \frac{\partial }{\partial t^{b}}:=\sum _{c}A_{ab}^{c}(t_{S})\frac{\partial }{\partial t^{c}}
\end{equation}

\subsection{Flat metrics}

Next ingredient needed for construction of the solution to WDVV-equation is
a flat metric on \( \mathcal{U} \) compatible with multiplication (\ref{circ}).
In this subsection we show that the period map \( \Psi ^{S} \) induces such
a metric under condition that the semi-infinite subspace \( S \) is isotropic
with respect to certain pairing.

Let us notice that sheaf \( \mathcal{R}_{\widehat{F}(t)} \) is equipped with
natural pairing between the opposite fibers: 
\[
G(\hbar ):\, \, H^{n}(X,\text {Re}\, \frac{\widehat{F}(t)}{\hbar }\rightarrow -\infty ;\mathbb {C})\otimes H^{n}(X,\text {Re}\, -\frac{\widehat{F}(t)}{\hbar }\rightarrow -\infty ;\mathbb {C})\rightarrow \mathbb {C}\]
 which is given by the Poincare pairing dual to intersection pairing on relative
cycles. Hence by continuity one has natural pairing on elements of \( L(t) \),
\( t\in \mathcal{U} \).

\begin{lem}
\label{G(L)}For \( \mu (\hbar ),\eta (\hbar )\in L(t)\subset H \) 
\[
G(\mu ,\eta )\in \hbar ^{n}\mathbb {C}[[\hbar ]]\]

\end{lem}
\begin{proof}
One has for \( \hbar  \) from any sector \( a<\arg (\hbar )<b \), \( b-a<\pi  \)
and \( \mu (\hbar ),\eta (\hbar )\in L(t) \): 
\[
G(\mu ,\eta )=\sum _{i,j}\#(\Delta _{i}^{+}\cap \Delta ^{-}_{j})(\int _{\Delta ^{+}_{i}}\mu (\hbar ))(\int _{\Delta ^{-}_{j}}\eta (-\hbar ))\]
 where sets \( \{\Delta _{i}^{+}\} \),\( \{\Delta _{j}^{-}\} \) are frames of relative
cycles in \( H_{n}(X,\text {Re}\, \frac{\widehat{F}(t)}{\hbar }\rightarrow -\infty ;\mathbb {C}) \),
\( H_{n}(X,\text {Re}\, -\frac{\widehat{F}(t)}{\hbar }\rightarrow -\infty ;\mathbb {C}) \)
correspondingly. Expansion for the integrals \( \int _{\Delta ^{+}_{i}}\mu
(\hbar )) \), \( \int _{\Delta ^{-}_{j}}\eta (-\hbar ) \)
at \( \hbar \rightarrow 0 \) can be evaluated using the steepest-descent method.
Recall that given a non-degenerate critical point \( p \) of \( \widehat{F}(t) \)
and a metric one has two relative cycles \( \Delta ^{+}_{p} \), \( \Delta ^{-}_{p} \)
from \( H_{n}(X,\text {Re}\, \frac{\widehat{F}(t)}{\hbar }\rightarrow -\infty ;\mathbb {C}) \)
and \( H_{n}(X,\text {Re}\, -\frac{\widehat{F}(t)}{\hbar }\rightarrow -\infty ;\mathbb {C}) \)
respectively which are formed by gradient lines \( \rho (r) \) of \( \text {Re}\, \widehat{F}(t)/\hbar  \)
such that \( \rho (r)\rightarrow p \) as \( r\rightarrow +\infty  \) and \( \rho (r)\rightarrow p \)
as \( r\rightarrow -\infty  \) respectively. Then one has the following asymptotic
expansions 
\begin{equation}
\label{int(delta+-)}
\int _{\Delta ^{+}_{p}}e^{\frac{\widehat{F}(x;t)}{\hbar }}(\varphi _{0}+\hbar \varphi _{1}+\ldots )=e^{\frac{\widehat{F}(p;t)}{\hbar }}(\text {const}^{(+)}\, \hbar ^{\frac{n}{2}}+O(\hbar ^{\frac{n}{2}+1}))
\end{equation}

\[
\int _{\Delta ^{-}_{p}}e^{-\frac{\widehat{F}(x;t)}{\hbar }}(\varphi _{0}-\hbar \varphi _{1}+\ldots )=e^{-\frac{\widehat{F}(p;t)}{\hbar }}(\text {const}^{(-)}\, \hbar ^{\frac{n}{2}}+O(\hbar ^{\frac{n}{2}+1}))\]
 Recall that \( f \) has \( (n+1) \) non-degenerate critical points \( p_{m} \),
\( p_{m} \): \( x_{0}=\ldots =x_{n}=\exp \frac{2\pi \sqrt{-1}m}{n+1} \) and
that the relative cycles \( \{\Delta ^{+}_{p_{m}}\} \), (resp. \( \{\Delta ^{-}_{p_{m}}\} \)
) form basis in \break \( H_{n}(X,\text {Re}\, \frac{f}{\hbar }\rightarrow -\infty ;\mathbb {C}) \),
(resp. \( H_{n}(X,\text {Re}\, -\frac{f}{\hbar }\rightarrow -\infty ;\mathbb {C}) \)).
It follows that 
\[
G(\mu ,\eta )=\text {const}\, \hbar ^{n}+O(\hbar ^{n+1})\]

\end{proof}
Let us assume now that \( S \) is isotropic with respect to \( G \) in the
following sense: 
\begin{equation}
\label{isotr}
\forall \, \, \, \mu (\hbar ),\eta (\hbar )\in S\, \, \, \, G(\mu ,\eta )\in \hbar ^{n-2}\mathbb {C}[[\hbar ^{-1}]]
\end{equation}
 An example of such subspace is considered in the next section. It is the subspace
generated by elements \( (\hbar \alpha )^{k}\hbar ^{-m-(n+1)\alpha } \), \( k\in \{0,\ldots ,n\} \),
\( m\in \mathbb {N} \).

The period map \( \Psi ^{S} \) induces a pairing on \( \mathcal{TM} \): 
\begin{equation}
\label{<>S}
\frac{\partial }{\partial t^{a}}\otimes \frac{\partial }{\partial t^{b}}\rightarrow G(\frac{\partial \Psi ^{S}}{\partial t^{a}},\frac{\partial \Psi ^{S}}{\partial t^{b}})
\end{equation}

\begin{prop}
\label{pr:g} For isotropic \( S \) the pairing induced by \( \Psi ^{S} \)
satisfies 
\begin{equation}
\label{<>h(n-2)}
G(\frac{\partial \Psi ^{S}}{\partial t^{a}},\frac{\partial \Psi ^{S}}{\partial t^{b}})=\hbar ^{n-2}g_{ab}
\end{equation}
 with \( g_{ab}(t) \) symmetric and non-degenerate. In the coordinates \( \{t^{a}_{S}\} \)
induced by the map \( \Psi ^{S}(t,\hbar =\infty ) \) the \( 2 \)-tensor \( g_{ab} \)
is constant. 
\end{prop}
\begin{proof}
The property (\ref{<>h(n-2)}) follows from \( \partial _{a}\Psi ^{S}(t,\hbar )\in S\cap \hbar ^{-1}L(t) \).
Notice that for any \( \mu (\hbar ),\eta (\hbar )\in H \): 
\[
G(\mu ,\eta )(\hbar )=(-1)^{n}G(\eta ,\mu )(-\hbar )\]
 It follows that \( g_{ab} \) is symmetric. Non-degeneracy of \( g_{ab} \)
follows from the analogous property of \( G \).

Notice that \( G \) is obviously \( \mathbb {C}[\hbar ^{-1}]- \)linear and
therefore one has induced pairing \( G_{0} \) with values in \( \mathbb {C}\cdot \hbar ^{n-2} \)
on elements of \( S/\hbar ^{-1}S \). It follows from (\ref{<>h(n-2)}) that
\[
G(\frac{\partial \Psi ^{S}}{\partial t^{a}_{S}},\frac{\partial \Psi ^{S})}{\partial t^{b}_{S}})=G_{0}([\frac{\partial \Psi ^{S}}{\partial t^{a}_{S}}],[\frac{\partial \Psi ^{S}}{\partial t^{b}_{S}}])\]
 By definition of the coordinates \( \{t^{a}_{S}\} \) one has \( [\frac{\partial \Psi ^{S}}{\partial t^{a}_{S}}]=constant \)
as an element of \( S/\hbar ^{-1}S \). 
\end{proof}
\begin{prop}
\label{gcirc} Metric \( g_{ab} \) is compatible with multiplication defined
by \( A_{ab}^{c}(t_{S}) \): 
\begin{equation}
\label{eq:gcirc}
\forall \, a,b,c\, \, \, \, \, \, g(a\circ b,c)=g(a,b\circ c)
\end{equation}

\end{prop}
\begin{proof}
Notice that \( G \) is locally constant with respect to the Gauss-Manin connection
\( \nabla ^{\mathcal{M}} \). Therefore 
\[
\partial _{b}G(\partial _{a}\Psi ,\partial _{c}\Psi )=G(\partial _{b}\partial _{a}\Psi ,\partial _{c}\Psi )+G(\partial _{a}\Psi ,\partial _{b}\partial _{c}\Psi )\]
 In the coordinates \( \{t^{a}_{S}\} \) one has 
\[
\partial _{b}G(\partial _{a}\Psi ,\partial _{c}\Psi )=0,\, \, \, \partial _{a}\partial _{b}\Psi =\hbar ^{-1}\sum _{d}A_{ab}^{d}(t_{S})\partial _{d}\Psi \]
 The pairing \( G \) satisfies \( G(\hbar ^{-1}\mu ,\eta )=-G(\mu ,\hbar ^{-1}\eta )=\hbar ^{-1}G(\mu ,\eta ) \)
\end{proof}
Let us denote 
\[
A_{abc}(t_{S}):=g(\partial _{a}\circ \partial _{b},\partial _{c})\]
 Corollary \ref{cor:AA,dA} and proposition \ref{gcirc} imply now that 
\[
A_{abc}(t_{S})=\partial _{a}\partial _{b}\partial _{c}\Phi ^{S}(t_{S})\]
 for some function \( \Phi ^{S} \) defined up to addition of terms of order
\( \leq 2 \), and that the function \( \Phi ^{S}(t_{S}) \) satisfy the system
of WDVV-equations: 
\begin{equation}
\forall \, \, a,b,c,d\, \, \, \, \, \sum _{e,f}\partial _{a}\partial _{b}\partial _{e}\Phi ^{S}g^{ef}\partial _{f}\partial _{c}\partial _{d}\Phi ^{S}=\sum _{e,f}\partial _{a}\partial _{d}\partial _{e}\Phi ^{S}g^{ef}\partial _{f}\partial _{b}\partial _{c}\Phi ^{S}
\end{equation}

\subsection{Symmetry vector field.}

Let us notice that the subspace \( L(t) \) satisfies 
\[
\mathcal{D}_{\partial /\partial \hbar }\, L(t)\subseteq \hbar ^{-2}L(t)\]
 Assume now that the subspace \( S\subset H \) satisfies the condition 
\begin{equation}
\label{Dstable}
\mathcal{D}_{\partial /\partial \hbar }\, S\subseteq \hbar ^{-1}S
\end{equation}
 Notice that this condition is enough to check only on finite set of elements
from \( S \) which form a basis of \( S/\hbar ^{-1}S \). Assume that one has
also \( \mathcal{D}_{\partial /\partial \hbar }\, \Omega _{0}\in \hbar ^{-1}S \). 

\begin{prop}
\label{pr:Euler} The periods \( \Psi ^{S,k}(t^{S}) \) satisfy 
\begin{equation}
\label{eq:Euler}
\frac{\partial \Psi ^{S,k}}{\partial \hbar }=\hbar ^{-1}E^{a}(t_{S})\frac{\partial \Psi ^{S,k}}{\partial t_{S}^{a}}
\end{equation}
 for some vector field \( E=\sum _{a}E^{a}(t^{S})\partial _{a} \)
\end{prop}
\begin{proof}
It follows from \( \frac{\partial \Psi ^{S,k}}{\partial \hbar }\in \hbar ^{-2}L(t)\cap \hbar ^{-1}S \).
The arguments are similar to the arguments from the proof of the proposition
\ref{Pr_pic-fuchs}. 
\end{proof}
\begin{prop}
The vector field \( E \) acts as a conformal symmetry on the multiplication
\( \, \circ \,  \) and the metric \( g_{ab} \), that is, for any vector fields
\( u,v\in \mathcal{T}_{\mathcal{M}} \)
\begin{equation}
\label{eq:Ecirc}
[E,u\circ v]-[E,u]\circ v-u\circ [E,v]=u\circ v
\end{equation}

\begin{equation}
\label{eq:Eg}
\text {Lie}_{E}g(u,v)-g([E,u],v)-g(u,[E,v])=(2-n)g(u,v)
\end{equation}

\end{prop}
\begin{proof}
It follows from \( [\nabla ^{\mathcal{M}},\mathcal{D}^{\mathcal{R}_{\widehat{F}(t)}}_{\partial /\partial \hbar }]=0 \)
and \( \mathcal{D}^{\mathcal{R}_{\widehat{F}(t)}}_{\partial /\partial \hbar }G=0 \)
respectively. 
\end{proof}

\subsection{Flat identity vector field\label{sident}. }

Assume that the subspace \( S \) satisfies the condition 
\begin{equation}
\label{omegaS}
\hbar ^{-1}\Omega _{0}\in S
\end{equation}
 Let \( e \) denotes vector field on \( \mathcal{M} \) which is affine in
the coordinates \( \{t_{S}\} \) and which corresponds to the class of \( \hbar ^{-1}\Omega _{0} \)
in the quotient \( S/\hbar ^{-1}S \). We can assume that \( e=\frac{\partial }{\partial t_{S}^{0}} \)
where \( t_{S}^{0} \) is one of the affine coordinates.

\begin{prop}
\label{pr:ident}The vector field \( e=\frac{\partial }{\partial t_{S}^{0}} \)
is the identity with respect to the multiplication \( \circ  \). 
\end{prop}
\begin{proof}
One has \( \hbar ^{-1}\Psi ^{S}\in \hbar ^{-1}L(t) \), \( \hbar ^{-1}\Psi ^{S}\in S \).
Therefore 
\[
\hbar ^{-1}\Psi ^{S}(t_{S})=\frac{\partial \Psi ^{S}(t_{S})}{\partial t_{S}^{0}}\]
 Differentiating this equality with respect to \( \frac{\partial }{\partial t_{S}^{a}} \)
and using eq.(\ref{Aabc}) we see that \( \frac{\partial }{\partial t_{S}^{a}}\circ \frac{\partial }{\partial t_{S}^{0}}=\frac{\partial }{\partial t_{S}^{a}} \). 
\end{proof}
In summary, we have

\begin{prop}
\label{pr:frob} For any semi-infinite subspace \( S \) satisfying the properties
(\ref{transv}), (\ref{isotr}), (\ref{Dstable}), (\ref{omegaS}) we have constructed
Frobenius manifold structure on \( \mathcal{U} \). 
\end{prop}
For the definition of Frobenius structure see \cite{D},\cite{KM}.

\begin{rem}
Similarly, there exists a family of solutions to WDVV-equations parameterized
by an open domain in isotropic affine grassmanian of semi-infinite subspaces
given an arbitrary similar pair \( (X,f) \) satisfying the following homological
Calabi-Yau type condition: there exists \( \Omega _{0}\in \Gamma _{Zar}(X,\Omega ^{n}) \),
\( n=\text {dim}_{\mathbb {C}}X \) such that the restriction on \( [e^{\frac{f}{\hbar }}\Omega _{0}] \)
of the symbol of Gauss-Manin connection induces isomorphism \( T_{[f]}\mathcal{M}\simeq (\hbar ^{-1}L([f]))/L([f]) \),
here under \( T_{[f]}\mathcal{M} \) we understand the corresponding space of
first-order deformations i.e. the vector space \( \Gamma _{Zar}(X,\mathcal{O}_{X})/\{Lie_{v}f\, |\, v\in \Gamma _{Zar}(X,\mathcal{T}_{X})\} \).
In particular this is true under some mild conditions whenever \( X \) admits
a nowhere vanishing ``holomorphic'' volume element. Also similar statement,
which deals with total sum of cohomology of complex of sheaves \( (\Lambda ^{*}TX_{Zar},\{f,\cdot \}) \)
and subspace of exponential integrals in total cohomology \( H^{*}(X,\text {Re}\, \frac{f}{\hbar }\rightarrow -\infty ;\mathbb {C}) \),
holds true. We assume in the previous statement unobstructedness of the corresponding
moduli space, which becomes non-trivial in such more general setting as deformations
become more of a noncommutative flavour. 
\end{rem}
\begin{rem}From some point of view the construction mentioned in the first part
of previous remark gives a certain global analogue of K. Saito and M.Saito theory
(\cite{SaitoKM}). This theory provides a construction of a Frobenius manifold
starting from a germ of holomorphic function having \emph{one} isolated singularity.
The construction uses microlocal algebraic analysis and, in particular, solution
to a certain Riemann-Hilbert problem obtained from some deep results from theory
of D-modules. On one hand in some cases, like quasi-homogenious singularity,
the local data involving germ of isolated singularity can be identified with
some global data associated with an algebraic function \( \mathbb {A}^{n}\rightarrow \mathbb {A}^{1} \)
of the type mentioned in the previous remark. It can be shown that any Frobenius
manifold arises from some abstract semi-infinite variation of Hodge structure
\( L(t) \) provided with some choice of constant opposite semi-infinite subspace
\( S \). Thus it is plausible that in such cases when global and local data
coincide and for specific choices of semi-infinite subspace \( S \) the construction
mentioned in the previous remark reproduces solutions to WDVV-equations which
were found in (\cite{SaitoKM}). Approach from previous remark can help to clarify
nature of solutions from (\cite{SaitoKM}). In a sense semi-infinite geometry
gives here an intrinsic tool for solution of a class of Riemann-Hilbert problems.
It should be noted also that solutions to WDVV-equations found in (\cite{SaitoKM})
are difficult to access for explicit calculations. Use of semi-infinite geometry
helps in principle to overcome this problem. On the other hand, applying
construction from the first part of the previus remark to general global data
 when, in particular, \( f \) has several singularities one can try to see
how the Saito's constructions associated with every singularity of \( f \)
get ``glued'' together somehow. Use of semi-infinite geometry could be helpful
in clarifying this gluing phenomena and finding global analogue \footnote{%
Based on example of quasi-homogenious polynomial existence of some global analogue
of Saito's theory was conjectured also in \cite{G}. Reader should be warned
that conjecture from \cite{G} must be corrected as the above examples show.
Namely, ``volume form'' \( \omega _{\lambda } \) (in notations of \cite{G})
with prescribed properties does not exist in general, instead there should exist
a family of top degree forms \( \omega _{\lambda }^{(0)}+\hbar \omega _{\lambda }^{(1)}+\ldots  \)
depending on the parameter \( \hbar  \). 
} of Saito's theory in most general setting of arbitrary pairs \( (X,f) \) \footnote{%
\emph{Note added:} After appearance of present paper C.Sabbah, to whom I would
like to express my gratitude for this interesting remark, told me about his
recent work (\cite{Sab}) containing extension of parts of Saito's theory in
certain global cases. However, global analogue of part of Saito's theory, which
involves proof of existence of ``primitive form'' in families depending on
parameters from moduli spaces like \( \mathcal{M}^{(X,F)} \) above and consequently
gives solutions to WDVV-equations, does not exist in literature to my knowledge. 
}, which could in turn provide a nice simplest example for study of general semi-infinite
variations of Hodge structures of Calabi-Yau type. \end{rem}

\section{Rational Gromov-Witten invariants of \protect\( \mathbb {C}\mathbb {P}^{n}\protect \)
and oscillating integrals. }

Consider the subspace \( S_{0}\subset H \) generated by elements 
\[
(\hbar \alpha )^{k}\hbar ^{-m-(n+1)\alpha },\, \, \, k\in \{0,\ldots ,n\},\, \, \, m\in \mathbb {N}\]
 Let us also put \( \Omega _{0}=[e^{\frac{f}{\hbar }}\frac{dx_{0\ldots }dx_{n}}{d(x_{0\ldots }x_{n})}]\in L_{0} \).
In this section we show that for \( S=S_{0} \) the quasi-homogeneous solution
to WDVV equation \( \partial ^{3}\Phi ^{S_{0}}(t_{S_{0}}) \) constructed in
the previous section coincides with the generating function of Gromov-Witten
invariants of \( \mathbb {C}\mathbb {P}^{n} \). The idea is to identify some
numerical invariants of these two solutions and then use reconstruction theorem
- I from \cite{KM}.

\begin{prop}
The subspace \( S_{0} \) has the properties (\ref{transv}), (\ref{isotr}),
(\ref{Dstable}), (\ref{omegaS}). 
\end{prop}
\begin{proof}
Properties (\ref{Dstable}) (\ref{omegaS}) are easy to check. In order to check
the property (\ref{transv}) let us notice that the subspace \( L(0) \) is
generated by elements 
\[
\hbar ^{2k+m}\frac{\partial ^{k}\xi }{\partial \hbar ^{k}}\, \text {mod}\, \alpha ^{n+1}\, \, \, k\in \{0,\ldots ,n\},\, \, m\in \mathbb {N}\cup \{0\}\]
 and that 
\[
\hbar ^{(n+1)\alpha }\frac{\partial ^{k}\xi }{\partial \hbar ^{k}}\, \text {mod}\, \alpha ^{n+1}=\text {const}\, \hbar ^{-k}\alpha ^{k}+O(\hbar ^{-(k+1)})\]
 To check the property (\ref{isotr}) let us notice that since \( \{\alpha ^{k}\} \)
is dual to locally constant basis \( \{\Delta _{k}\} \) and \( \alpha ^{i}\hbar ^{-(n+1)\alpha } \)
are single-valued sections therefore for any \( i,j \) the pairing \( G\bigl (\alpha ^{i}\hbar ^{-(n+1)\alpha },\alpha ^{j}\hbar ^{-(n+1)\alpha }\bigr ) \)
does not depend on~\( \hbar  \): 
\begin{equation}
\notag \forall \, \, \, \, i,j\, \, \, \, G\bigl (\alpha ^{i}\exp (-(n+1)\log (\hbar )\alpha ),\alpha ^{j}\exp (-(n+1)\log (\hbar )\alpha )\bigr )=\text {const}
\end{equation}
 Differentiating this equality with respect to \( \frac{\partial }{\partial \hbar } \)
one gets 
\[
G\bigl (\alpha ^{i+1}\hbar ^{-(n+1)\alpha },\alpha ^{j}\hbar ^{-(n+1)\alpha }\bigr )=G\bigl (\alpha ^{i}\hbar ^{-(n+1)\alpha },\alpha ^{j+1}\hbar ^{-(n+1)\alpha }\bigr )\]
 Therefore \( G\bigl (\alpha ^{i}\hbar ^{-(n+1)\alpha },\alpha ^{j}\hbar ^{-(n+1)\alpha }\bigr )=0 \)
for \( i+j>n \). 
\end{proof}
Let us consider system of local coordinates \( \{t^{k}_{S_{0}}\} \) , \( k=0\ldots n \),
induced by the map (\ref{Psi(infty)}) from linear coordinates on \( S/\hbar ^{-1}S \)
corresponding to the basis given by classes of elements \( \alpha ^{k}\hbar ^{k-1-(n+1)\alpha } \)
modulo \( \hbar ^{-1}S \). Notice that this notation is compatible with notation
for the flat identity from \S \ref{sident} since \( [\hbar ^{-1}\Omega _{0}]=[\hbar ^{-1-(n+1)\alpha }] \).
We put below \( y^{k}=t^{k}_{S_{0}} \)to simplify notations.

\begin{lem}
\label{E0} In coordinates \( y^{k} \) the symmetry vector field is written
as 
\[
E(y^{k})=(\sum _{k=0}^{n}(k-1)y^{k}\frac{\partial }{\partial y^{k}})-(n+1)\frac{\partial }{\partial y^{1}}\]

\end{lem}
\begin{proof}
One has 
\[
\Psi ^{S_{0}}(y)=\Omega _{0}+\sum _{k=0}^{n}\hbar ^{k-1-(n+1)\alpha }\alpha ^{k}y^{k}\, \, \, mod\, \, \, \hbar ^{-1}S_{0}\]
 Therefore \begin{multline}
\frac{\partial\Psi^{S_{0}}(y)}{\partial\hbar}
=\\=-(n+1)\alpha\hbar^{-1-(n+1)\alpha} +
\sum_{k=0}^{n}\hbar^{k-2-(n+1)\alpha}\alpha^{k} (k-1)y^{k}
\,\,\,mod\,\,\,\hbar ^{-2}S_{0}\notag
\end{multline}
\begin{equation}
\notag \frac{\partial \Psi ^{S_{0}}(y)}{\partial y^{k}}=\hbar ^{k-1-(n+1)\alpha }\alpha ^{k}\, \, \, mod\, \, \, \hbar ^{-1}S_{0}
\end{equation}
 and 
\[
\frac{\partial \Psi ^{S_{0}}(y)}{\partial \hbar }=\hbar ^{-1}\bigl (-(n+1)\frac{\partial \Psi ^{S_{0}}(y)}{\partial y^{1}}+\sum _{k=0}^{n}(k-1)y^{k}\frac{\partial \Psi ^{S_{0}}(y)}{\partial y^{k}}\bigr )\, \, \, mod\, \, \, \hbar ^{-2}S_{0}\]
 This is an exact equality by proposition \ref{pr:Euler}
\end{proof}
\begin{lem}
\label{pr:t=3D3D3D3D3D0} 
\[
\frac{\partial \Psi ^{S_{0},i}(y)}{\partial y^{k}}|_{y=0}=\frac{1}{\hbar (n+1)^{k}}\int _{\Delta _{i}}f^{k}\exp (\frac{1}{\hbar }f)\frac{dx_{0\ldots }dx_{n}}{d(x_{0\ldots }x_{n})}\, \, mod\, \, \hbar L(0)\]
 In other words at \( y=0 \) the elements \( \frac{\partial }{\partial y^{k}} \)
correspond to elements \( (\frac{\sum ^{n}_{i=0}x_{i}}{n+1})^{k}\in \Gamma _{Zar}(X,\mathcal{O}_{X})/\{Lie_{v}(f)\, |\, v\in \Gamma _{Zar}(X,\mathcal{T}_{X})\} \)
from \( T_{[f]}\mathcal{M}^{(X,f)} \)
\end{lem}
\begin{proof}
One has 
\[
\hbar ^{k-1}(\hbar \frac{\partial }{\partial \hbar })^{k}\Psi ^{S_{0}}|_{y=0}=\hbar ^{k-1}(-(n+1)\alpha )^{k}\hbar ^{-(n+1)\alpha }\, \, mod\, \, \hbar ^{-1}S_{0}=\]

\[
=(-(n+1))^{k}\frac{\partial \Psi ^{S_{0}}(y)}{\partial y^{k}}|_{y=0}\, \, mod\, \, \hbar ^{-1}S_{0}\]
 we see that 
\[
\hbar ^{k-1}(\hbar \frac{\partial }{\partial \hbar })^{k}\Psi ^{S_{0}}|_{y=0}-(-(n+1))^{k}\frac{\partial \Psi ^{S_{0}}(y)}{\partial y^{k}}|_{y=0}\in \hbar ^{-1}L(0)\cap \hbar ^{-1}S_{0}\]
 Therefore \( \hbar ^{k-1}(\hbar \frac{\partial }{\partial \hbar })^{k}\Psi ^{S_{0}}|_{y=0}=(-(n+1))^{k}\frac{\partial \Psi ^{S_{0}}(y)}{\partial y^{k}}|_{y=0} \)
and 
\[
\frac{\partial \Psi ^{S_{0},i}(y)}{\partial y^{k}}|_{y=0}=\frac{1}{\hbar (n+1)^{k}}\int _{\Delta _{i}}f^{k}\exp (\frac{1}{\hbar }f)\frac{dx_{0\ldots }dx_{n}}{d(x_{0\ldots }x_{n})}\, \, mod\, \, \hbar L(0)\]

\end{proof}
\begin{cor}
\label{A(t=3D3D3D3D3D3D3D0)}At \( y=0 \) the multiplication (\ref{circ})
is given by 
\[
A_{ij}^{k}(0)=\delta _{i+j-k\, mod\, n+1,0}\]
 in the basis \( \{\frac{\partial }{\partial y^{k}}\} \). 
\end{cor}
\begin{proof}
\( \Gamma _{Zar}(X,\mathcal{O}_{X})/\{Lie_{v}(f)\, |\, v\in \Gamma _{Zar}(X,\mathcal{T}_{X})\}\simeq \mathbb {C}[p]/(p^{n+1}-1) \)
with \( p=[\frac{\sum ^{n}_{i=0}x_{i}}{n+1}] \). 
\end{proof}
Let us calculate now the pairing (\ref{<>S}) written in coordinates \( y^{k} \).
By proposition \ref{pr:g} and lemma \ref{pr:t=3D3D3D3D3D0} the value of \( \hbar ^{-n+2}G(\frac{\partial \Psi ^{S_{0}}}{\partial y^{i}},\frac{\partial \Psi ^{S_{0}}}{\partial y^{j}}) \)
is equal to the coefficient in front of \( \hbar ^{n} \) in \begin{multline} G([\frac{1}{(n+1)^{i}}f^{i}\exp
(\frac{1}{\hbar }f)\frac{dx_{0}\ldots dx_{n}}{d(x_{0}\ldots
x_{n})}],[\frac{1}{(n+1)^{j}}f^{j}\exp (\frac{1}{\hbar
}f)\frac{dx_{0}\ldots dx_{n}}{d(x_{0}\ldots x_{n})}])=\\=
 \sum_{
 m=0}^{n}\frac{1}{(n+1)^{i+j}}(\int_{\Delta_{p_{m}}^{+}} f^{i}
 e^{\frac{f}{\hbar}}\frac{dx_{0}\ldots dx_{n}}{d(x_{0}\ldots
x_{n})}))(\int_{\Delta_{p_{m}}^{-}} f^{j}
 e^{-\frac{f}{\hbar}}\frac{dx_{0}\ldots dx_{n}}{d(x_{0}\ldots
x_{n})}))\label{sum(pm)}
\end{multline} where \( \{p_{m}\} \) is the set of critical points of \( f \),
\( p_{m}:x_{0}=\ldots =x_{n}=\zeta ^{m}, \)\( \zeta =\exp \frac{2\pi \sqrt{-1}}{n+1} \);
and \( \Delta _{p_{m}}^{+},\Delta _{p_{m}}^{-} \) are the associated relative
cycles defined above in section 3. The lowest order coefficient in (\ref{sum(pm)})
can be calculated by expanding \( f \) near critical points up to second order:
\begin{multline} \frac{1}{(n+1)^{i+j}}\sum_{
 m=0}^{n} \int_{\Delta_{p_{m}}^{+}} f^{i}(p_{m})
 \exp\bigl(\frac{1}{\hbar}(f(p_{m})+\zeta^{m}\sum_{1 \leq k,l\leq n}(1-\frac{x_{k}}{\zeta^{m}})
 (1-\frac{x_{l}}{\zeta^{m}}))\bigr)
 \frac{dx_{1}}{x_{1}}\ldots\frac{dx_{n}}{x_{n}}
])\\ \int_{\Delta_{p_{m}}^{-}} f^{j}(p_{m})
 \exp\bigl(\frac{1}{\hbar}(f(p_{m})+\zeta^{m}\sum_{1 \leq k,l\leq n}(1-\frac{x_{k}}{\zeta^{m}})
 (1-\frac{x_{l}}{\zeta^{m}}))\bigr)
 \frac{dx_{1}}{x_{1}}\ldots\frac{dx_{n}}{x_{n}}
])\sim\\ \sim\hbar ^{n}(\sum_{
 m=0}^{n}
 \zeta^{m(i+j-n)}) + O(\hbar^{n+1})\sim\hbar ^{n}\delta
_{i+j-n,0}+ O(\hbar^{n+1})
\end{multline}

\begin{cor}
\label{g0} The pairing on \( T\mathcal{U} \) given by 
\[
<\frac{\partial }{\partial t^{i}_{S_{0}}},\frac{\partial }{\partial t^{j}_{S_{0}}}>=\delta _{i+j,n}\]
 satisfies the conditions (\ref{eq:gcirc}) and (\ref{eq:Eg}). 
\end{cor}
Let us consider Taylor expansion of the constructed solution \( \Phi ^{S_{0}}(y) \)
at \( y=0 \), where \( \partial _{ijk}^{3}\Phi ^{S_{0}}(y)=\sum _{m}\delta _{i+m,n}A_{jk}^{m}(y) \)
and we assume that the Taylor expansion of \( \Phi ^{S_{0}}(y) \) has no terms
of order \( \leq 2 \).

\begin{lem}
\label{expans} The Taylor expansion of \( \Phi ^{S_{0}} \) has the form 
\begin{equation}
\label{tayl}
\Phi ^{S_{0}}=\frac{1}{2}(\sum _{ij}y^{i}y^{j}y^{0}\delta _{i+j,n})+\sum _{m_{a}}\frac{\sigma (m_{1},\ldots ,m_{n})}{m_{1}!\ldots m_{n}!}(y^{1})^{m_{1}}\ldots (y^{n})^{m_{n}}
\end{equation}
 with 
\begin{equation}
\label{Etaylor}
\sigma (m_{1}+1,m_{2},\ldots ,m_{n})=\frac{3-n+\sum _{k=2}^{n}(k-1)m_{k}}{n+1}\sigma (m_{1},m_{2},\ldots ,m_{n})
\end{equation}
 and if \( \frac{3-n+\sum _{k=2}^{n}(k-1)m_{k}}{n+1} \) is not integral then
\begin{equation}
\label{sigma0}
\sigma (m_{1},m_{2},\ldots ,m_{n})=0
\end{equation}

\end{lem}
\begin{proof}
The form of the term which contains \( y^{0} \) follows from lemma \ref{pr:t=3D3D3D3D3D0}
and proposition \ref{pr:ident}. Equation (\ref{Etaylor}) follows from (\ref{eq:Ecirc}),(\ref{eq:Eg}).
Therefore expansion (\ref{tayl}) can be rewritten as expansion in powers of
\( (y^{0},e^{\frac{y^{1}}{n+1}}=1+\frac{y^{1}}{n+1}+\ldots ,y^{2},\ldots ,y^{n}) \).
To prove that this expansion contains only integral powers of \( e^{y^{1}} \)
let us introduce an additional parameter. Namely, let us consider family \( (X_{q},f) \)
where \( f=x_{0}+\ldots +x_{n} \) and \( X_{q}\subset \mathbb {C}^{n} \) is
defined by equations \( x_{0}\cdot \ldots \cdot x_{n}=q \), \( q\in \mathbb {C} \),
\( |q-1|<1 \). Identifying \( X_{q} \) with \( X \) via multiplication of
every coordinate by \( q^{-\frac{1}{n+1}} \) we see that \( (X_{q},f)\sim (X,q^{\frac{1}{n+1}}f) \)
as deformations of \( (X,f) \). It follows from (\ref{L0}) that 
\[
\sum ^{n}_{k=0}\alpha ^{k}\int _{\Delta _{k}}\exp (q^{\frac{1}{n+1}}\frac{f}{\hbar })\frac{dx_{0\ldots }dx_{n}}{d(x_{0\ldots }x_{n})}]=\hbar ^{-(n+1)\alpha }q^{\alpha }(1+O(\frac{1}{\hbar }))\]
 we see that \( [\exp (q^{\frac{1}{n+1}}\frac{f}{\hbar })\frac{dx_{0\ldots }dx_{n}}{d(x_{0\ldots }x_{n})}]\in S_{0}+[\Omega _{0}] \)
and that 
\begin{equation}
\label{OmegaQ}
[\exp (q^{\frac{1}{n+1}}\frac{f}{\hbar })\frac{dx_{0\ldots }dx_{n}}{d(x_{0\ldots }x_{n})}]-[\Omega _{0}]=\ln q\, \alpha \hbar ^{-(n+1)\alpha }\, \text {mod}\, \hbar ^{-1}S
\end{equation}
 Therefore \( (X_{q},f) \) corresponds to \( (0,y_{1},0,\ldots ,0)\in \mathcal{U} \),
\( q=e^{y^{1}} \), in the moduli space of deformations of \( (X,f) \). One
can now repeat the above story starting with \( (X_{q},f) \), \( q\in \mathbb {C}\setminus \{0\} \)
and consider element \( \Psi ^{S_{0}}(y,\hbar ;q) \) representing periods of
the variation of semi-infinite Hodge structures associated with moduli space
of deformations of \( (X_{q},f) \). The normalization condition is given by
\[
\Psi ^{S_{0}}(y,\hbar ;q)-[\exp (q^{\frac{1}{n+1}}\frac{f}{\hbar })\frac{dx_{0\ldots }dx_{n}}{d(x_{0\ldots }x_{n})}]\in S_{0}\]
 and the parameters \( \{y^{i}\} \) are specified via 
\[
\Psi ^{S_{0}}(y,\hbar ;q)=[\exp (q^{\frac{1}{n+1}}\frac{f}{\hbar })\frac{dx_{0\ldots }dx_{n}}{d(x_{0\ldots }x_{n})}]+\sum _{k=0}^{n}\alpha ^{k}\hbar ^{k-1-(n+1)\alpha }y^{k}\, \, \, mod\, \, \, \hbar ^{-1}S_{0}\]
 Notice that this is correctly defined since \( \exp (2\pi \sqrt{-1}\alpha )(S_{0})=S_{0} \)
and therefore \( S_{0} \) is invariant under monodromy around \( q=0 \). It
follows from (\ref{OmegaQ}) that 
\[
\Psi ^{S_{0}}(y_{0},y_{1},\ldots ,y_{n},\hbar ;q)=\Psi ^{S_{0}}(y_{0},0,\ldots ,y_{n},\hbar ;e^{y^{1}}q)\]
 Therefore expansion of \( \partial ^{3}\Phi ^{S_{0}}(y) \) can contain only
integral powers of \( e^{y^{1}} \)which gives (\ref{sigma0}). 
\end{proof}
\begin{thm}
\label{th:main}The Taylor expansion (\ref{tayl}) of the function \( \Phi ^{S_{0}}(y) \)
coincides with generating series for Gromov-Witten invariants of \( \mathbb {C}\mathbb {P}^{n} \):
\begin{equation}
\label{main}
\Phi ^{S_{0}}(y)=\mathcal{F}^{\mathbb {C}\mathbb {P}^{n}}(y)
\end{equation}

\end{thm}
\begin{proof}
It follows from lemmas \ref{E0}, \ref{expans}, corollaries \ref{A(t=3D3D3D3D3D3D3D0)},
\ref{g0} and theorem 3.1 from \cite{KM}
\begin{rem}
(This remark is a part of a joint work with M.Kontsevich) It is an interesting
open question whether a homological mirror symmetry conjecture can be formulated
in this setting in order to explain the new mirror phenomena. A natural analog
of derived category of coherent sheaves associated with \( (X,f) \) does not
seem to catch all the necessary information since it is most likely that it
has no non-trivial objects. A resolution of this problem might be in introducing
a further generalization of \( A_{\infty }- \)categories where objects exist
only in some virtual sense. On the other hand, homological mirror symmetry conjecture
in opposite direction which identifies \( D^{b}Coh(\mathbb {P}^{n}) \) with
a version of derived Fukaya category can be formulated. A sketch of the appropriate
definition of the derived Fukaya category associated with \( (X,f) \) can be
found in \cite{Knotes} (some technical details of the construction were checked
in recent preprints math.SG0007115 \& 0010032 by P.Seidel). Whether there exists
some \emph{numerical} mirror symmetry related with the opposite homological
mirror symmetry conjecture (i.e. \( D^{b}Coh(\mathbb {P}^{n}) \) versus a Fukaya
category) is an interesting open question deserving further study. 
\end{rem}
\end{proof}


\begin{thebibliography}{EHX}
\bibitem[B1]{B1}S.~Barannikov, \emph{Generalized periods and mirror symmetry in
higher dimensions}, Preprint IHES/M/99/16,
math.AG/9903124; \emph{Extended moduli spaces and mirror symmetry in dimensions
$N>3$.} Ph.D. thesis, UC Berkeley, Spring 1999.
\bibitem[B2]{B2}S.~Barannikov, \emph{Quantum periods --I. Semi-infinite variations of Hodge
structures,} preprint ENS DMA-00-19, alg-geom/0006193, to appear in IMRN.
\bibitem[BK]{BK}S.\ Barannikov, M.\ Kontsevich, \emph{Frobenius manifolds and formality of Lie
algebras of polyvector fields}, Internat.\ Math.\ Res.\ Notices, no.4 (1998),
201-215. 
\bibitem[D]{D}B.~Dubrovin, \emph{Geometry of 2d topological field theories}, LNM 1620, Springer,
1996, 120-348. 
\bibitem[EHX]{EHX}T.Eguchi,K.Hori, Ch.-Sh. Xiong, \emph{Gravitational quantum cohomology,} Int.J.Mod.Phys.
A12 (1997), 1743-1782. 
\bibitem[GKZ]{Kapr}I.~Gelfand, M.~Kapranov, A.Zelevinsky, \emph{Hypergeometric functions, toric
geometry and Newton polyhedra}. Special functions (Okayama 1990), 104-121, Springer,
1991. 
\bibitem[G]{G}A.Givental, \emph{Homological geometry and mirror symmetry.} In Proc. of ICM-94
Zurich, Birkhauzer, 1995, 472-480.
\bibitem[K1]{K1}M.~Kontsevich, \emph{Deformation quantization of Poisson manifolds--I}, preprint
q-alg/9709040. 
\bibitem[K2]{Knotes}M.~Kontsevich, \emph{Lectures at ENS (Paris), Spring 1998,} (notes taken by
J.Bellaiche, J.-F.Dat, I.Marin, G.Racinet and H.Randriambololona). 
\bibitem[KM]{KM}M. Kontsevich, Yu. Manin, \emph{Gromov-Witten classes, quantum cohomology, and
enumerative geometry}, Commun. Math. Phys. \textbf{164} (1994), 525-562. 
\bibitem[M]{M}Yu.I.~Manin, \emph{Sixth Painlev\'{e} equation, universal elliptic curve and
mirror of \( P^{2} \),} Amer. Math. Soc. Transl. ser.2,186, 131-151. 
\bibitem[PS]{PS}A.Pressley, G.Sigal, \emph{Loop groups,} Clarendon Press, Oxford, 1988. 
\bibitem[Sab]{Sab}C.~Sabbah, \emph{Hypergeometric periods for a tame polynomial.} C.R. Acad. Sci.
Ser. I Math., 328 (1999), 603-608.
\bibitem[S,S]{SaitoKM}K.~Saito, \emph{Period mapping associated to a primitive form.} Publ. RIMS,
Kyoto Univ. 19 (1983), 1231-1264. M.~Saito, \emph{On the structure of Brieskorn
lattices}. Ann Inst. Fourier, 39 (1989), 27-72. 
\bibitem[SW]{SW}G.Sigal, G.Wilson, \emph{Loop groups and equations of KdV type,}-Publ.Math.
IHES, 61, 1985, 5-65
\end{thebibliography}
\end{document}